\documentclass[a4paper,12pt]{amsart}

\usepackage{amsxtra}
\usepackage{amssymb}
\usepackage[all]{xypic}
\usepackage{epsfig}

\xyoption{dvips}

\theoremstyle{definition}
\newtheorem{ntn}{Notation}

\newtheorem{rem}[ntn]{Remark}
\newtheorem{exa}[ntn]{Example}
\theoremstyle{plain}

\newtheorem{prp}[ntn]{Proposition}
\newtheorem{thm}[ntn]{Theorem}

\theoremstyle{remark}

\DeclareMathAlphabet{\mathds}{U}{dsrom}{m}{n}
\DeclareMathAlphabet{\mathsc}{U}{rsfs}{m}{n}

\DeclareMathOperator{\Gal}{Gal}
\DeclareMathOperator{\Spec}{Spec}

\newcommand{\Q}{\mathbb{Q}}
\newcommand{\N}{\mathbb{N}}
\newcommand{\Z}{\mathbb{Z}}
\newcommand{\F}{\mathbb{F}}
\renewcommand{\P}{\mathbb{P}}
\newcommand{\R}{\mathbb{R}}
\newcommand{\A}{\mathbb{A}}
\newcommand{\C}{\mathbb{C}}

\renewcommand{\O}{\mathcal O}

\renewcommand{\a}{\alpha}

\newcommand{\lr}{\rightarrow}

\begin{document}

\title{Fields of Definition of Rational Points on Varieties}
\author{Jordan Rizov}
\address{Mathematisch Instituut\\ P.O. Box 80.010\\ 3508 TA Utrecht\\ The Netherlands}
\email{rizov@math.uu.nl}
\thanks{I thank Ben Moonen, Frans Oort, Andrea Giacobbe, Ivan Chipchakov and Grigor Grigorov for stimulating discussions.}
\maketitle

\begin{abstract}
Let $X$ be a scheme over a field $K$ and let $M_X$ be the intersection of all subfields $L$ of $\bar K$ such that $X$ has a $L$-valued point. In this note we
prove that for a variety $X$ over a field $K$ finitely generated over its prime field one has that $M_X = K$.
\end{abstract}
\vspace{23pt}
Let $K$ be a field and fix an algebraic closure $K \subset \bar K$. For a scheme $X$ over $K$ denote by $\mathsc C_X$ the collection of all fields 
$K \subset L \subset \bar K$ such that $X$ has a $L$-valued point $x\colon \Spec(L) \lr X$. Define the field $M_X$ as
$$
 M_X = \bigcap_{L \in \mathsc C_X} L
$$
where the intersection takes place in $\bar K$. We will
be interested in how big $M_X$ can be and in particular in some cases in which it is $K$
itself.

If the set of $K$-rational points $X(K)$ is non-empty, then obviously $M_X$ is $K$. In general, if $K$ is a perfect field and $X$ is a scheme of finite type over $K$ then $M_X$ is a finite 
Galois extension of $K$. Indeed, let $x \colon \Spec(L) \lr X$ be a $L$-valued point on $X$ corresponding to a point $x$ on the topological space $X$ and an
inclusion $\kappa(x) \lr L$ where $\kappa(x) = \O_x/{\mathfrak m}_x$ is the residue field of $x$ (see \cite[Ch. II, \S 2, Exercise 2.7]{HAG} and \cite[Ch. II,
\S 4, Prop. 3]{Mum-RedBk}). For any $\sigma \in \Gal(\bar \Q/K)$ one has
the conjugate $x^\sigma$ of $x$ which is a $\sigma(L)$-valued point on $X$. Hence for any field $L \in \mathsc C_X$ and any $\sigma \in \Gal(\bar \Q/K)$ the field
$\sigma(L)$ is also in $\mathsc C_X$ which implies that $M_X/K$ is a Galois extension.
  
\begin{rem}\label{observation}
Suppose given two schemes $X_1$ and $X_2$ defined over $K$ and a morphism $f \colon
X_1 \lr X_2$ over $K$. Then one has that $\mathsc C_{X_1} \subset \mathsc C_{X_2}$ and therefore $M_{X_2} \subset M_{X_1}$. In particular if $M_{X_1} = K$, then the field $M_{X_2}$ is $K$, as well.
\end{rem}

Before going on let us consider some illuminating examples under different conditions
imposed on $X$ and $K$.
\begin{exa} 
If we take $X$ to be $\Spec(K)$, then clearly $M_X = K$. 
\end{exa}
\begin{exa}\label{dobarpr}
Consider the curve $X \colon x^2 + y^2 + z^2 = 0$ in $\P^2$ over $\Q$ (i.e. $K =
\Q$). One can take $P_1 = [i:0:1]$ which has field of definition $\Q(i)$ and the point $P_2 =
[\sqrt{-2}:1 :1]$  giving the field $\Q(\sqrt{-2})$. Clearly, the intersection
of those two fields is $\Q$, so $M_X = \Q$. This is an example where $X$ is a
non-singular, projective curve of
genus 0 over $\Q$ with no $\Q$-rational point.
\end{exa}
\begin{exa}\label{dobarkp}
Let $K = \F_q$ be a finite field and let $X$ be a non-singular, quasi-projective curve over $K$. We may assume that $X$ is contained in its complete,
non-singular model $X'$ over $K$. Let $X = X'\setminus \{P_1,\dots,P_r\}$ (the complete case is treated in the same way) and denote the genus of $X'$ by
$g$. If $n \in \N$ and $N_{q^n}$ denote the number of $\F_{q^n}$ rational points on $X'$ then by the Weil bound we have that
$$
 N_{q^n} \geq 1+ q^n - 2g\sqrt{q^n}.
$$
Therefore, if $n$ is sufficiently large one has that $N_{q^n} \geq r + 1$ and hence $X(\F_{q^n})$ is not empty. Choose two natural numbers $n_1$ and $n_2$ which are
sufficiently large so that $X(\F_{q^{n_i}})$ is not empty for $i = 1,2$ and $\gcd(n_1,n_2) = 1$. Then we have that
$$
 M_X \subset \F_{q^{n_1}} \cap \F_{q^{n_2}} = \F_q
$$
and hence $M_X = K$. 
\end{exa}
\begin{exa}\label{loshpr}
Consider again the curve $X \colon x^2 + y^2 + z^2 = 0$ in $\P^2$ but take this time $K$ to
be $\R$. Then we have that $M_X$ is $\C$ since $X$ has no $\R$-valued points.
\end{exa}
\begin{exa}\label{patologichen}
Take $K=\Q$ and consider the polynomial $f(x) = x^3 - 7x +7$. It is irreducible
over $\Q$ since it has no rational zeros. Let $\a = \a_1,
\a_2$ and $\a_3$ be its roots. Since the discriminant of $f(x)$ is $7^2$ its
Galois group is isomorphic to $\Z/3\Z$ and the field $M = \Q(\a)$ is a Galois extension of
$\Q$ of degree 3. Let $P_i = (\a_i,\a_i^2)$ for $i=1,2,3$ be three points in $\A^2_{\bar \Q}$
and consider the three lines passing through them:
$$
\begin{matrix}
l_1 = y - (\a_1 + \a_2)x + \a_1\a_2 = 0 & {\rm passing\ through}\ P_1\ {\rm
and}\ P_2; \\
l_2 = y - (\a_1 + \a_3)x + \a_1\a_3 = 0 & {\rm passing\ through}\ P_1\ {\rm
and}\ P_3; \\
l_3 = y - (\a_2 + \a_3)x + \a_2\a_3 = 0 & {\rm passing\ through}\ P_2\ {\rm
and}\ P_3.
\end{matrix}
$$
Define the scheme $X \subset \A^2_{\bar \Q}$ to be given by the equation $l_1l_2l_3 = 0$.
The Galois group $\Gal(\bar \Q/ \Q)$ permutes the three points $P_1, P_2$ and
$P_3$ and respectively the three lines $l_1, l_2$ and $l_3$ in $\A_{\bar \Q}^2$. Hence $X$ is defined over $\Q$ and it is irreducible over $\Q$. 

Let $P \colon \Spec(L) \lr X$ be a $L$-valued point on $X$ for some field $L
\subset \bar \Q$. If $\sigma \in \Gal(\bar \Q/L)$, then it fixes the point $P$ on
$X_{\bar \Q}$. Since an automorphism in $\Gal(\bar \Q/\Q)$, acting on $\A_{\bar \Q}^2$, permutes the three lines we
see that, acting on $X_{\bar \Q}$, it has fixed points if and only if it acts as the trivial
permutation on $\{l_1,l_2,l_3\}$. Hence $\sigma$ must leave the points $P_i$,
$i=1,2,3$ fixed and therefore $M$, as well. Thus we conclude that $M \subset L$. Since
$X$ has $M$-valued points (the points $P_i$ for $i=1,2,3$) we have that $M_X = M = \Q(\a)$.

In this example one can take $f(x)$ to be any irreducible polynomial over $\Q$ with Galois group isomorphic to $\Z/3\Z$.
\end{exa}
The last two examples suggest that if the field of definition $K$ is `too big' or the
scheme $X$ is somehow `too bad' then the field $M_X$ is a non-trivial extension of $K$. On the other
hand the argument given in example \ref{dobarpr} can be easily generalized
to number fields and other curves.

Before stating the main result let us make the following convention: In this note a variety $X$ over a field $K$ will mean a separated, geometrically integral scheme
$X$ of finite type over $K$. In particular, $X$ is geometrically irreducible. Also from now on we will assume that $K$ is a finitely generated field over its
prime field. We know that for those fields Faltings' finiteness theorem holds (see \cite[Ch. I, \S 2]{L-NT3}). Further, if char($K$)$=0$ or if char($K$)$=p$ and ${\rm tr.deg}_{\F_p}K
\geq 1$, then Hilbert's irreducibility theorem holds for the field $K$. We refer to \cite[Ch. 9]{Lang-DG} and more precisely to Theorem 4.2 and the remark following
it.
\begin{thm}\label{stranna}
Let $K$ be a finitely generated field over its prime field and let $X$ be a variety over $K$. Then one has that $M_X = K$.
\end{thm}
\begin{proof}
{\bf Step 1.} 
We will first show that it is enough to consider non-singular, quasi-projective varieties. If $X$ is not complete, then by Nagata's 
compactification theorem one can find a complete variety $\bar X$ and an open immersion $i \colon X \hookrightarrow \bar X$. By Chow's Lemma there 
exits a projective variety $Y'$ over $K$ and a birational isomorphism $\pi' \colon Y' \lr \bar X$. Let $Y$ be an alteration of $Y'$ (see \cite[\S 1 and \S 4,
Thm. 4.1]{deJ-Alt}), let 
$\pi \colon Y \lr \bar X$ be the composition morphism and let $X' = \pi^{-1}(i(X))$. Then by Remark \ref{observation} we have that $M_X \subset M_{X'}$.
Hence it is enough to show the validity of the theorem assuming that $X$ is a non-singular, quasi-projective variety over $K$. 
\newline
\newline 
We may assume that $X \subset \P^N$ for some $N$. In the next two steps we will show here that it is enough to 
prove the theorem assuming that $\dim X =1$.
\newline
\newline
{\bf Step 2.}
Suppose that $K$ is an infinite field. If $\dim X = 1$ then the result follows from Proposition \ref{rezkrivi} below. Suppose that $\dim X = m \geq 2$. Then by Bertini's Theorem (\cite[Ch. II, \S 8, Thm. 8.18)]{HAG}) we know that the set $U$ of points $u$ in the dual
projective space $\check \P^N$ corresponding to hyperplanes $H \subset \P^N_{\kappa(u)}$ such that $H\cap X$ is
smooth of dimension $m-1$ over the residue field $\kappa(u)$ of $u$ contains a dense open subset of $\check \P^N$. Since $K$ is infinite the intersection $U \cap
\check \P^N(K)$ is non-empty. Hence one can find a hyperplane $H$ defined over $K$ satisfying Bertini's Theorem. Further, by \cite[Ch. III, \S 11, Exercise 11.3]{HAG} the
intersection $H\cap X$ is geometrically connected and hence it is geometrically irreducible or in other words it is a quasi-projective variety of
dimension $m-1$ over $K$. Repeating this $\dim X - 1$ times one can find a non-singular, quasi-projective curve $Y \subset X$ defined over $K$. By 
Remark \ref{observation} one has that $M_X \subset M_Y$. Now the claim follows from Proposition \ref{rezkrivi} below.
\newline
\newline
{\bf Step 3.}
Let $K$ be a finite field. The case $\dim X = 1$ was considered in example \ref{dobarkp}. Assume that $\dim X \geq 2$. We will find again a quasi-projective
curve defined over $K$ contained in $X$ by intersecting $X$ with hypersurfaces. By Theorem 3.3 and Remarks (1) and (2) in \cite{P-BFF} there exists a geometrically
integral, smooth hypersyrface $H \subset \P^N$ defined over $K$ such that the intersection $X \cap H$ is a smooth variety of dimension $\dim X - 1$. Repeating
this $\dim X - 1$ times we can find a non-singular, quasi-projective curve $Y$ in $X$ defined over $K$. Then just like in Step 2 we conclude the claim from
Remark \ref{observation} and example \ref{dobarkp}.
\end{proof}
\begin{prp}\label{rezkrivi}
Let $K$ be an infinite field which is finitely generated over its prime field. If $X$ is a quasi-projective curve defined over $K$, then $M_X = K$.
\end{prp}
\begin{proof}
We will split up the proof into three steps.
\newline
\newline
{\bf Step 1.} Assume that $X$ is a complete, non-singular curve of genus at least 2 and there is a morphism $f \colon X \lr \P^1$ over $K$ of prime degree $p$.
Hilbert's irreducibility theorem assures that there are infinitely many points
$P \in \P^1(K)$ such that the fiber $f^{-1}(P) = \{Q_1, \dots , Q_r\}$ consists of points which
are defined over extensions $K(Q_i)$ of $K$ of degree $p$. If among all
those fields (for all points $P\in \P^1(K)$ as above), there are two which are different, then their intersection will be
$K$ (as they do not have non-trivial subfields). Hence we would have that
$M_X=K$. Assume that all fields $K(Q_i)$ for all $P \in \P^1(K)$ as before are
the same. Then we have infinitely many
points on $X$ defined over a fixed extension $L = K(Q_i)$ of $K$. As $X$ is of genus at least 2 we get a contradiction with
Faltings' finiteness theorem. Thus we conclude that
$M_X$ is $K$ in this case.
\newline
\newline
{\bf Step 2.} Now assume that $X$ is complete and non-singular. In general, one should not expect to be able to find a morphism as in Step 1. Instead, we will
construct a covering $\pi \colon X' \lr X$ over $K$ for some curve $X'$ satisfying the
assumptions of Step 1. Then we could conclude the claim of the proposition using Remark \ref{observation}. Such a curve
can be viewed as a divisor on $X \times \P^1$ so we will look at special divisors on this ruled surface.

Let $a$ be a natural number which we will fix later and consider the divisor $D(a) = 2X + a\P^1$ on
$X \times \P^1$. Following the notations of \cite[Ch V, \S 2]{HAG} we put $(X,X)_{X\times \P^1} = -e$. Then using Proposition
2.3, Lemma 2.10 and Corollary 2.11 of \cite[Ch. 5, \S 2]{HAG}, one sees that
$$
 (D(a), X) = (2X + a\P^1, X) = a - 2e
$$
and the `adjunction formula' for the divisor $D(a)$ has the form
$$
 (D(a), D(a) + K_{X \times \P^1}) = 2a + 2(2g_X - 2 - e)
$$
where $K_{X \times \P^1}$ is the canonical class of $X \times
\P^1$. Let us choose $a$ so that
$$
\begin{matrix}
a > 2e & \\
a - 2e & {\rm is\ a\ prime\ number} \\
a + (2g_X - 2 - e) \geq 1. &
\end{matrix}
$$
The first condition ensures that the linear system $|D(a)|$ contains a non-singular, geometrically irreducible curve $X'$ defined over $K$. Indeed, one uses \cite[Ch.
V, \S 5, Cor. 2.18]{HAG}. Hartshorne assumes that $K$ is algebraically closed. Since the proof only deals with very ample line bundles and uses Bertini's Theorem 
\cite[Ch. II, \S 8, Thm. 8.18]{HAG} and \cite[Ch. III, \S 11, Exercise 11.3]{HAG} it remains valid over $K$, as $K$ is an infinite field.

The degree of the morphism $f \colon X' \hookrightarrow X \times \P^1 \lr
\P^1$ is exactly $(X',X) = a - 2e$ which is a prime number and by the adjunction formula the genus of $X'$ is $a + (2g_X -
2 - e) + 1$ which is at least $2$ by our choice. Hence $X'$ satisfies the conditions of Step 1. Therefore by Remark \ref{observation} we conclude that $M_X =
M_{X'} = K$.
\newline
\newline
{\bf Step 3.} Let $X$ be as in the proposition. If it is not complete one can find a completion $X'$ of $X$ which
is also defined over $K$. Take the normalization $X''$ of $X'$ and let $\tilde X \subset X''$ be the preimage of $X$. The curve $X''$ is non-singular, 
projective and defined over $K$. One can now apply Step 1 and 2 above to
$X''$.  Clearly those proofs, and more precisely the one of Step 1, can be carried
over excluding a finite number of points (which of course should not change the
field of definition $K$). In other words one sees that $M_{\tilde X} = K$. Hence by Remark \ref{observation} we have that 
$M_X \subset M_{\tilde X} = K$ and therefore $M_X=K$.
\end{proof}

\begin{rem}
Note that thought we distinguished the two cases $K$ is finite and $K$ is infinite the two proofs go exactly in the same lines. One tries to find
finite extensions $L_1$ and $L_2$ of $K$ which are `different' as subfields of $\bar K$, such that $X(L_i)$ is not empty for $i= 1,2$ and so that one can control the intermediate
fields $K \subset M \subset L_i$. In the case $K$ is finite this is easily achievable using the Weil bound. If $K$ is infinite one makes use of Hilbert's
irreducibility theorem instead. Below we shall present a proof based on a completely different idea. Namely, in the case $K$ is a number field one tries to find
sufficiently many prime ideals of $K$ splitting completely in $M_X$. This proof was suggested to us by Grigor Grigorov. 
\end{rem}

Let K be a number field and let $\O_K$ be the ring of integers in $K$. For a prime ideal $\mathfrak p$ of $K$ denote by $k_\mathfrak p$ the 
residue field $\O_K/\mathfrak p$. Let $q_\mathfrak p$ be the number of
elements in $k_\mathfrak p$. Then by definition one has that the norm $N(\mathfrak p)$ of $\mathfrak p$ is $q_\mathfrak p$. Denote by $K_\mathfrak p$ the completion of $K$ at 
$\mathfrak p$ and let $\O_{K_\mathfrak p}$ be ring of integers in $K_\mathfrak p$.
\newline
\newline
{\it Proof of Proposition} \ref{rezkrivi} {\it assuming that $K$ is a number field.}   
We already saw that one can assume that $X$ is non-singular and it is contained in its complete non-singular model $X'$ defined over $K$. 
We have that $X = X' \setminus \{P_1,\dots,P_m\}$ for some $m \in \N$ (the proof in the complete case is the same). Take a projective embedding of $X'$ over $K$ in to $\P^N_K$ for some $N$ and 
its flat closure $\mathcal X'$ over $\O_K$ in $\P^N_{\O_K}$. Let $\mathcal X$ be the complement $\mathcal X'\setminus \{\mathcal P_1,\dots,\mathcal P_m\}$
where $\mathcal P_i$, $i = 1,\dots,m$, is the flat closure of $P_i$ over $\O_K$. Then there is a finite set of primes $\Sigma$ such that $\mathcal X'$ is smooth over 
$U = \Spec(\O_K) \setminus \Sigma$. For a prime ideal $\mathfrak p \not\in \Sigma$ let $N_\mathfrak p$ be the number of points in $\mathcal X'(k_\mathfrak p)$. 
The Weil bound reads 
$$
N_\mathfrak p \geq 1+q_\mathfrak p - 2g \sqrt{q_\mathfrak p}.
$$
where $g$ is the genus of $X'$. Hence if $q_\mathfrak p = N(\mathfrak p)$ is sufficiently large one has that $N_\mathfrak p \geq m+1$. So enlarging $\Sigma$, if needed, we may assume that 
$\mathcal X(k_\mathfrak p)$ is not empty for all $\mathfrak p \not \in \Sigma$.

Fix a prime ideal $\mathfrak p \not\in\Sigma$. Since $\mathcal X$ is smooth over $U$ and $\mathcal X(k_\mathfrak p)$ is non-empty one can apply Hensel's lemma
(see \cite[\S 2.3, Prop. 5]{NM}) to 
conclude that $X(K_\mathfrak p)$ is non-empty. Therefore by Theorem 1.3 in \cite{MB-Sk2} one can find a finite extension $L$ of $K$ such that
$\mathfrak p$ splits completely in $L$ and $X$ has a $L$-valued point. Hence $\mathfrak p$ splits completely in $M_X$. Thus all but finitely many ideals 
(at most those in $\Sigma$) split completely in $M_X$. By Corollary 6.6 in \cite[Ch. V, \S 6]{N-CFT} we have that $M_X = K$.
\qed
\begin{rem}
Theorem \ref{stranna} could be viewed as a variant of Theorem 5.1 in \cite[\S 5]{D-ShV1} where Deligne proves that for a Shimura datum $(G,X)$ and any finite
extension $L$ of the reflex field $E(G,X)$ of the Shimura variety $Sh(G,X)$ there exists a special point $x$ on $X$ such that its reflex field $E(x)$ is
linearly disjoint from $L$. This result is used in proving the uniqueness of the canonical model of $Sh(G,X)$ over $E(G,X)$. We came across the main result of this note
considering a similar descent problem. 
\end{rem}

How big can $M_X$ be in general? We already saw in examples
\ref{loshpr} and \ref{patologichen} that depending on $X$ and $K$
the field $M_X$ can be a non-trivial extension of $K$. Using the
construction in example \ref{patologichen} one can find $X$ over
$\Q$ such that $[M_X:\Q]$ is arbitrary large. On the other hand if
$X$ is a non-singular, projective curve defined over a field $K$, then $l(K_X) = g$, where $g$ is the
genus of $X$ and $K_X$ is its canonical class. If $g \geq 2$ then
there is a non-constant $K$-rational function $f$ in $L(K_X)$. It
gives a morphism $f \colon X \lr \P^1$ of degree at most $\deg K_X
= 2g-2$. Hence there is a $L$-valued point for some extension
$L/K$ with $[L:K] \leq 2g-2$. Therefore we have that $[M_X:K] \leq
2g-2$.

\bibliographystyle{amsalpha}
\providecommand{\bysame}{\leavevmode\hbox to3em{\hrulefill}\thinspace}
\providecommand{\MR}{\relax\ifhmode\unskip\space\fi MR }
\providecommand{\MRhref}[2]{%
  \href{http://www.ams.org/mathscinet-getitem?mr=#1}{#2}
}
\providecommand{\href}[2]{#2}

\end{document}